%% file: oneface.tex
\documentclass [11pt,oneside]{article}

\reversemarginpar \textwidth=14cm \textheight=19cm

\usepackage{amssymb} \usepackage{amsfonts} \usepackage{amsmath}
\usepackage{amsthm} \usepackage{epsfig}

\newcommand{\mettifig}[1]{\epsfig{file=#1}}

\include{defs_oneface}

\author{Roberto \titsc{Frigerio} \and Bruno \titsc{Martelli}
\and Carlo \titsc{Petronio}}

\title{Complexity and Heegaard genus of\\ an infinite class of compact 3-manifolds}

\begin{document}

\maketitle

\begin{abstract}
	\noindent
	Using the theory of hyperbolic manifolds with totally geodesic 
	boundary, we provide
	for every $n\geqslant 2$ a class $\calM_n$ 
	of such manifolds all having Matveev
	complexity equal to $n$ and Heegaard genus equal to $n+1$. All the elements 
	of $\calM_n$ have a single boundary
	component of genus $n$, and $\#\calM_n$ grows at least exponentially with $n$.

  \vspace{4pt}

\noindent MSC (2000): 57M27 (primary), 57M20, 57M50 (secondary).
\end{abstract}

\noindent
This paper is devoted to the investigation of the class $\calM_n$ of 
orientable compact $3$-manifolds having an ideal
triangulation with $n\geqslant 2$ tetrahedra and a single edge. 
We show in particular for each $M$ in $\calM_n$ that 
the Heegaard genus of $M$ is equal to $n+1$, and that
the complexity of $M$ in the sense of Matveev
is equal to $n$.  Moreover we prove that the classical invariants,
such as homology and the Turaev-Viro invariants, cannot
distinguish two distinct members of $\calM_n$ from each other.
However, using the fact that each $M$ in $\calM_n$ 
carries a hyperbolic metric with totally geodesic boundary, we prove
that $M$ has a unique ideal triangulation with $n$ tetrahedra.
We exploit this property showing that the number of elements
of $\calM_n$ grows at least exponentially with $n$.  This implies in particular
the previously unknown fact that the rate of growth of the number of
orientable boundary-irreducible acylindrical manifolds of complexity $n$ is also
at least exponential in $n$.
The class $\calM_n$ was already considered in~\cite{Hea}, but none
of our results was covered there.

\section{Manifolds with a one-edged triangulation}
In this section we introduce the class of manifolds we are interested in, and
we prove their many remarkable topological and geometric properties.

\begin{figure}
\begin{center}
\mettifig{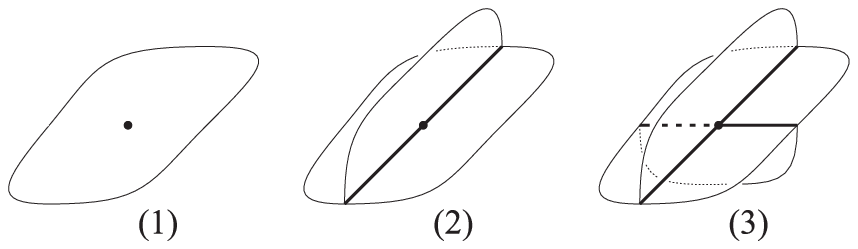,width=9cm}
\nota{Local aspect of a quasi-standard polyhedron.} \label{standard_nhbds:fig}
\end{center}
\end{figure}

\begin{figure}
\begin{center}
\mettifig{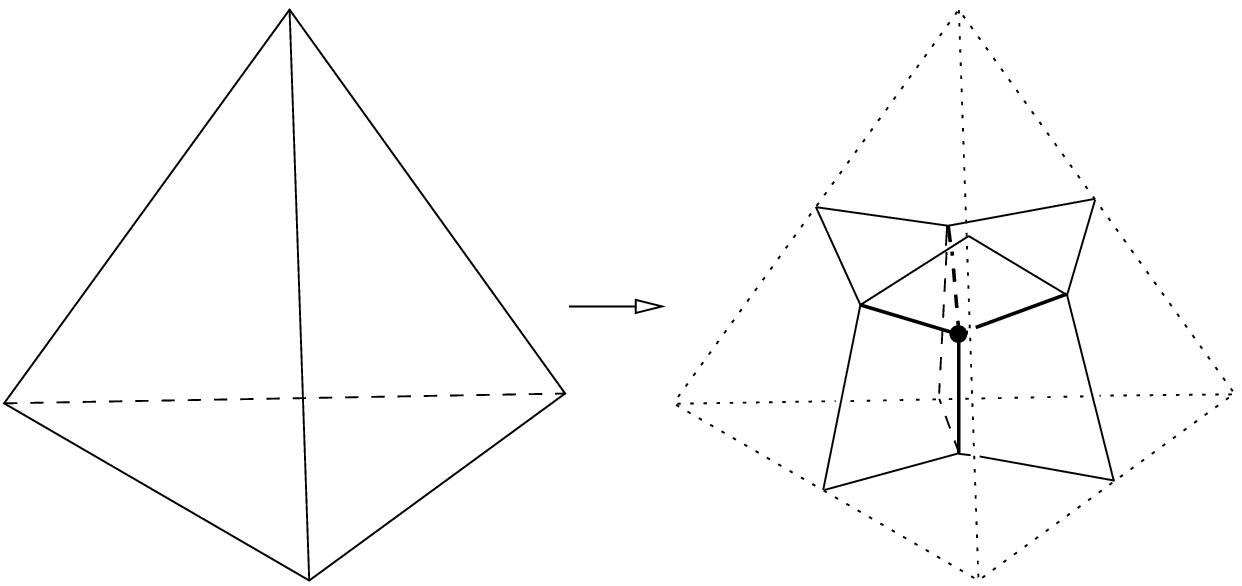, width=8 cm}
\nota{Duality between ideal triangulations and spines.} \label{dualspine:fig}
\end{center}
\end{figure}

\paragraph{Ideal triangulations and spines}
We begin by recalling some definitions. An \emph{ideal tetrahedron} 
is a tetrahedron with its vertices
removed. An \emph{ideal triangulation} of a compact $3$-manifold $M$ with boundary is a
realization of the interior of $M$ as a gluing of some ideal tetrahedra,
induced by a simplicial pairing of the faces. A \emph{spine} of $M$ is a compact
polyhedron $P$ such that $M\setminus P = \partial M\times [0,1)$. A 2-dimensional
polyhedron $Q$ is \emph{quasi-standard} if every point has a neighbourhood 
homeomorphic to one of the
polyhedra shown in Fig.~\ref{standard_nhbds:fig}.
We denote by $V(Q)$ the set of points 
having regular neighbourhoods of type (3),
and by $S(Q)$ the set of points 
having regular neighbourhoods of type (2) or (3). If $Q\setminus S(Q)$ consists of open
cells and $S(Q)\setminus V(Q)$ consists of open edges, we say that $Q$ is
\emph{standard}, and we call \emph{faces} the components of $Q\setminus S(Q)$.
An ideal triangulation of $M$ defines in a natural
way a dual standard polyhedron, which is in fact a spine of $M$ 
(see Fig.~\ref{dualspine:fig}).
We now define the class of manifolds investigated in this paper. 
For every integer $n\geqslant 2$
we set
\begin{eqnarray*}
\calM_n = \{M & : & \textrm{dim}(M)=3,\ M\textrm{ is compact and orientable},\ \partial M\neq\emptyset,\\
&& M \textrm{ admits an ideal triangulation with one edge and } n \textrm{ tetrahedra}\}.
\end{eqnarray*}
The class $\calM_n$ can be defined in various equivalent ways, as the next lemma shows.

\begin{lemma} \label{definitions:lem} Let $M$ be a compact $3$-manifold and let 
$\calT$ be an ideal triangulation of $M$ consisting of $n$ tetrahedra.
The following facts are equivalent:
\begin{enumerate}
\item $\calT$ has one edge;
\item $\chi(M) = 1-n$;
\item $\chi(M) = 1-n$ and $\partial M$ is connected.
\end{enumerate}
\end{lemma}
\begin{proof}
Set $\widehat M = M/\partial M$, let $|\partial M|$ be the number of components of $\partial M$,
and let $x$ be the number of edges of $\calT$.
We have $\chi(\widehat M) = \chi(M)-\chi(\partial M) + |\partial M|$ and
$\chi(\partial M) = 2\chi(M)$. Extending $\calT$ to a cellularization of $\widehat M$,
we also get

$$\chi(\widehat M) = |\partial M| - x + 2n - n = |\partial M| - x + n.$$
Summing up, we get $\chi(M) = x-n$, which shows $(1)\Leftrightarrow(2)$.

We are left to show that if $x=1$ then $|\partial M|=1$. Since there is a single edge,
we have $|\partial M|\leqslant 2$. Suppose $A$ and $B$ are distinct components of $\partial M$,
and examine a triangular face $F$ of $\calT$. Let 
$e_1$, $e_2$, and $e_3$ be the edges of $F$ viewed abstractly (\emph{i.e.}~before
the embedding in $M$). Since $e_1,e_2,e_3$ become the same edge in $M$, they should
all join $A$ to $B$, which is clearly impossible:
if $e_1$ and $e_2$
join $A$ to $B$, then $e_3$ joins either $A$ or $B$ to itself.
\end{proof}

\paragraph{Topological and geometric properties}
Before proving our main theorem, we recall the definition of Matveev complexity
of a compact 3-manifold with boundary, and the notion of hyperbolic 3-manifold with
geodesic boundary. A compact 2-dimensional
polyhedron $Q$ is said to be \emph{simple} if the link of every point in $Q$ is contained in
the 1-skeleton $K$ of the tetrahedron. (Note that a standard polyhedron is obviously simple.)
A point having the whole of $K$ as a link is 
called a \emph{vertex}, and its regular neighbourhood
is as shown in Fig.~\ref{standard_nhbds:fig}-(3). This implies that the set $V(Q)$ of the
vertices of $Q$ consists of isolated points, so it is finite.
The \emph{complexity} $c(M)$ of
a compact 3-manifold $M$ with boundary is the minimal number of vertices of
a simple spine of $M$.

A \emph{hyperbolic $3$-manifold with geodesic boundary} is a complete
Riemannian manifold with boundary
which is locally isometric to a half-space of hyperbolic $3$-space
$\matH^3$. By Mostow's
Rigidity Theorem (see~\cite{FriPe} for an explicit statement in the geodesic boundary case), 
every compact 3-manifold admits at most one hyperbolic structure
with geodesic boundary, and has therefore a well-defined hyperbolic volume
(if any). A useful tool for the computation of hyperbolic volumes, used in the sequel, is the
Lobachevsky function $L:\matR\to\matR$ defined by
$$ L(\omega)=-\int\limits_0^\omega \log |2 \sin u |\,\textrm{d}u.$$

To state our result we also recall that for any integer $r\geqslant 2$,
after fixing $q_0$ in $\matC$ such that $q_0^2$ is a primitive $r$-th root 
of unity, 
a real-valued invariant $\TV_r$
for compact 3-manifolds with boundary was defined by Turaev and Viro in~\cite{TuVi}.
We consider here these invariants normalized so that $\TV_r(S^3)=1$ for any $r$.
Their computation involves the complex-valued quantum $6j$-symbols
$\left\{\begin{array}{ccc} i & j & k \\ l & m & n \end{array}\right\}$, 
where $i,j,k,l,m,n$ are half-integers, and
the real-valued quantum integers $[k] = \frac{q_0^k-q_0^{-k}}{q_0-q_0^{-1}}$, defined for any
integer $k$.

\begin{teo}\label{main:teo}
Let $M\in \calM_n$. Then:
\begin{enumerate}
\item
$M$ is hyperbolic with geodesic boundary and its volume is given by
$${\rm vol}(M)= n\cdot\left( 8L \left(\frac{\pi}4\right) - 3\int\limits_0^{\frac{\pi}{3n}}
\mathrm{arccosh} \left(\frac{\cos t}{2\cos t -1}\right) \,{\rm d}t\right); $$
\item
$M$ is boundary-irreducible and acylindrical;
\item Every closed incompressible surface in $M$ is parallel to the boundary;
\item
$H_1(M;\matZ)\cong\matZ^n$;
\item
The Heegaard genus of $M$ is equal to $n+1$;
\item
$c(M)=n$;
\item
The $r$-th Turaev-Viro invariant of $M$ is given by
$$\TV_r(M) = \sum_{h\in\matN,\ 0\leqslant 3h\leqslant r-2} 
\left\{\begin{array}{ccc} h & h & h \\ h & h & h \end{array}\right\}^n\cdot
[2h+1]^{1-n}$$
\end{enumerate}
\end{teo}

\begin{figure}
\begin{center}
\mettifig{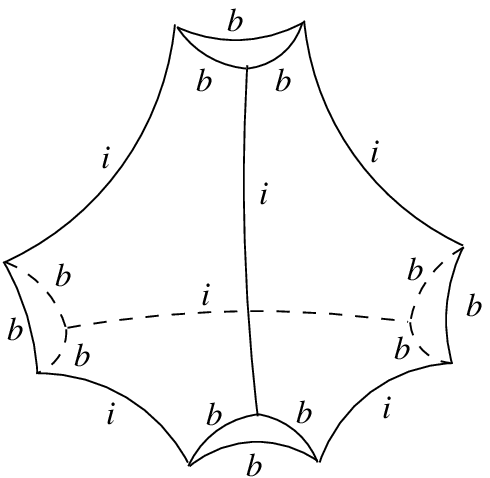, width=4cm}
\nota{Length of the edges of a regular truncated tetrahedron.} \label{regtrunc:fig}
\end{center}
\end{figure}

Before giving the proof of Theorem \ref{main:teo}, we introduce the notion of
(hyperbolic) \emph{truncated tetrahedron}~\cite{Fuj, Ko, FriPe}.
Let $\Delta$ be
a tetrahedron and let $\Delta^*$ be the combinatorial polyhedron obtained
by removing from $\Delta$ small open stars of the vertices. We call
\emph{lateral hexagon} and \emph{truncation triangle} the intersection of
$\Delta^*$ respectively with a face and with the link of a vertex
of $\Delta$.
The edges of the truncation triangles are called \emph{boundary edges},
the other edges of $\Delta^*$ are called \emph{internal edges}.
A \emph{hyperbolic truncated tetrahedron} 
is a realization of $\Delta^*$ as
a compact polyhedron in $\matH^3$,
such that the truncation triangles are geodesic triangles,
the lateral hexagons are geodesic hexagons, and truncation triangles and lateral
hexagons lie at right angles to each other. A truncated tetrahedron is
\emph{regular} if all the dihedral angles along its internal edges are
equal to each other. It turns out~\cite{Fuj,FriPe} that for every $\theta$ with
$0<\theta<\pi/3$ there exists up to isometry exactly one regular truncated
tetrahedron $\Delta^*_\theta$ of dihedral angle $\theta$. The boundary edges
of $\Delta^*_\theta$ all have the same length $b=b(\theta)$, 
and the internal edges all have the same length $i=i(\theta)$,
as shown in Fig.~\ref{regtrunc:fig}.
The volume of $\Delta^*_\theta$ is given by (see~\cite{Miy}):
$$ {\rm vol}(\Delta^*_\theta) = 8L \left(\frac{\pi}4\right) - 3\int\limits_0^{\theta}
\mathrm{arccosh} \left(\frac{\cos t}{2\cos t -1}\right) \,\textrm{d}t. $$
As the proof of Theorem \ref{main:teo} will now make clear,
truncated tetrahedra can be used as building blocks
to construct hyperbolic manifolds with geodesic 
boundary~\cite{Fuj, Ko,FriPe}.

\dimo{main:teo}
First of all, we fix 
a one-edged ideal triangulation $\calT$ of $M$
and we denote by $e$ the single edge of $\calT$. 

In order to give $M$ a hyperbolic structure, we identify each
tetrahedron of $\calT$ with a copy of the regular truncated
tetrahedron $\Delta^*_{\pi/3n}$. Due to the symmetries of
$\Delta^*_{\pi/3n}$, every pairing between the faces of the tetrahedra
of $\calT$ 
can be realized by an isometry between the corresponding lateral hexagons
of the $\Delta^*_{\pi/3n}$'s. 
This procedure defines a hyperbolic metric
on $M$, possibly with a cone singularity along $e$.
But each tetrahedron is incident 6 times to $e$, 
so the cone angle around $e$ is
$6\cdot n\cdot\pi/3n = 2\pi$, and $M$ is actually hyperbolic
without singularities. Finally,
${\rm vol}(M) = n\cdot {\rm vol}(\Delta^*_{\pi/3n})$ and point (1) is proved.

Point (2) is a general consequence of the existence
of a hyperbolic structure with geodesic boundary, and point (3)
is proved using Haken's theory of normal surfaces. If $\Sigma\subset M$
is a closed incompressible surface, then $\Sigma$ can be isotoped into 
normal position
with respect to $\calT$, so the intersection of $\Sigma$ with an ideal 
tetrahedron $T\in\calT$
consists of triangles and squares. Since $\Sigma$ intersects all the edges 
of $T$ in the
same number $k$ of points, it easily follows that $k$ is even and 
$\Sigma\cap T$
consists of $k/2$ parallel triangles at each vertex of $T$. Therefore 
$\Sigma$
consists of $k/2$ disjoint surfaces parallel to $\partial M$.

Concerning point (4),
let $P$ be the standard spine of $M$ dual to $\calT$. Since $M$ collapses onto $P$,
we have $H_1(M;\matZ)\cong H_1(P;\matZ)$, and
we can use cellular homology  to compute $H_1(P;\matZ)$.
To do so, we choose a maximal tree $Y$
in the 4-valent graph $S(P)$. Then $S(P)\setminus Y$ consists of $n+1$ edges
$e_1,\ldots,e_{n+1}$. Choose an orientation on each $e_i$ and on the single face
$F$ of $P$. The face $F$ is incident three times to each $e_i$. 
For $i=1,\ldots,n+1$ let $r_i$ be the sum of a contribution $\pm 1$
over the three instances $F$ runs along $e_i$, with sign depending on whether 
orientations are matched or not. So $r_i\in\{-3,-1,1,3\}$ and
$H_1(P;\matZ)\cong\matZ^{n+1}/\langle r\rangle$ where $r=(r_1,\ldots,r_{n+1})$. 
We will now prove
that $r_i=\pm 1$ for some $i$, which implies that $H_1(P;\matZ)\cong\matZ^n$.
Let $v\in V(P)$ be an extremal vertex of $Y$, \emph{i.e.}~a vertex adjacent to only one edge
of $Y$, whence to three of the $e_i$'s, say $e_{i_1}, e_{i_2}, e_{i_3}$
(indices could be repeated according to multiplicity of incidence).
Looking at the picture of the neighbourhood of $v$ in $P$,
shown in Fig.~\ref{standard_nhbds:fig}-(3), one now easily sees that $F$
cannot be incident three times in the same direction to each $e_{i_1}, e_{i_2}, e_{i_3}$, 
so we have $r_i=\pm 1$ for some $i\in\{i_1,i_2,i_3\}$.

Let us turn to point (5).  Lemma~\ref{definitions:lem} shows that
$\partial M$ has genus $n$ but, by point (2), $M$ is not a
handlebody, so its genus is at least $n+1$. A Heegaard surface having 
genus
$n+1$ is simply given by the boundary of a regular neighbourhood of
$\partial M \cup e$, whence the conclusion.

We now prove point (6).
The standard spine $P$ dual to $\calT$ has $n$ vertices, so we have $c(M)\leqslant n$.
By point (2), a result of Matveev~\cite{Mat}
shows that there is a standard spine $Q$ of $M$ (not just a simple one)
with precisely $c(M)$ vertices. An Euler characteristic computation gives
$\chi(Q)=x-c(M)$ where $x\geqslant 1$ is the number of faces of $Q$,
so $c(M)\geqslant 1-\chi(Q) = 1-\chi(M) = n$, the last equality having been proved in
Lemma~\ref{definitions:lem}.

Point (7) is an easy calculation. We follow the notation of~\cite{TuVi}:
since there is only one face $F$, a colouring of $P$ is given by assigning
to $F$ a half-integer in $\{0,1/2,1,3/2,\ldots,(r-2)/2\}$. Since the same $F$ 
is incident three
times to each edge, the colouring is admissible if it is an integer
$h$ with $0\leqslant 3h \leqslant r-2$. Each such colouring contributes
to $\TV_r(M)$ with a summand given by the product of a factor
$\left\{\begin{array}{ccc} h & h & h \\ h & h & h \end{array}\right\}\cdot[2h+1]^{-1}$
for each vertex and of a factor $[2h+1]$ due to the single face of $P$.
\finedimo

\begin{rem}\label{decomposition:rem}
{\em The proof of Theorem~\ref{main:teo}-(1) actually shows that
every one-edged ideal triangulation of $M\in\calM_n$ is combinatorially equivalent
to a decomposition of $M$ into $n$ regular truncated
tetrahedra of dihedral angle $\pi/3n$.}
\end{rem}

\begin{rem}\label{volume:rem}
{\em It was proved in~\cite{Miy} that among hyperbolic 3-manifolds 
with geodesic boundary and fixed 
Euler characteristic $\chi<0$, those having 
minimal volume are decomposed into $1-\chi$ copies of
$\Delta^*_{\pi/3(1-\chi)}$. Therefore $\calM_n$ is precisely the set of 
hyperbolic 3-manifolds $M$ with geodesic boundary
having minimal volume among orientable manifolds with $\chi(M) = 1-n$.}
\end{rem}

\begin{rem}\label{complexity:rem}
{\em It follows from the proof of Theorem~\ref{main:teo}-(6) that
$\calM_n$ is also the set of all hyperbolic manifolds $M$ having minimal
complexity among orientable manifolds with $\chi(M) = 1-n$.}
\end{rem}

\begin{rem}\label{TV:rem}
{\em  Point (7) of Theorem~\ref{main:teo} shows in particular that
$\TV_r(M)$, for $M$ in $\calM_n$, depends only on $n$ and $r$ (actually, on
$q_0$), but not on $M$. This could also be proved using the fact that
two standard spines
with the same incidence relations between faces and vertices produce the same
Turaev-Viro invariants (see~\cite{Mat:TV} for more details).
In fact, each manifold in $\calM_n$ admits a spine with one face and
$n$ vertices, so the incidence relations are always the same.}
\end{rem}

\begin{rem}\label{ori:rem}
{\em We believe that Theorem~\ref{main:teo} extends, with minor variations,
to non-orientable manifolds admitting a one-edged triangulation. However,
to prove this extension, one should first generalize Matveev's theorem~\cite{Mat}
on standardness of minimal spines to the 
case of non-orientable, boundary-irreducible manifolds not containing
projective planes or essential M\"obius strips.}
\end{rem}

\paragraph{Uniqueness of the minimal spine}
We are now left with two tasks: we must provide examples of
manifolds in $\calM_n$, and we must be able 
to distinguish manifolds in the same $\calM_n$. We will face the former task
in the next
section, by constructing standard polyhedra with $n$ vertices and one face. 
Concerning the latter task, we have just shown that
homology, Heegaard genus, Turaev-Viro invariants, and volume 
of manifolds in $\calM_n$ depend on $n$ only, so we need
a more powerful tool. This tool is provided by hyperbolic geometry. We recall
that the \emph{cut-locus} of a hyperbolic manifold $M$ with geodesic boundary
is the set of all points of $M$ which admit at least two distinct distance-minimizing
geodesics to $\partial M$.

\begin{teo}\label{distinct:teo}
Every $M\in\calM_n$ has a unique standard spine with $n$ vertices, homeomorphic
to the cut-locus of $M$.
\end{teo}
\begin{proof}
Let $C\subset M$ be the cut-locus of $M$, and let $P$ be a standard spine of $M$
with $n$ vertices. By Remark~\ref{decomposition:rem},
$P$ is dual to a decomposition $\calT$
of $M$ into $n$ regular truncated tetrahedra.
We claim that $C$ intersects each tetrahedron $T\in\calT$ 
as in Fig.~\ref{dualspine:fig}-right.
This implies that $C$ is homeomorphic to $P$, whence the conclusion.

To prove our claim it is sufficient to show that for every
tetrahedron $T$ of $\calT$, every point $p$ in $T$, and 
every distance-minimizing geodesic $\gamma$ connecting $p$ to
$\partial M$, we have that $\gamma$ is entirely contained in $T$.
If this were not true, since $\gamma$ meets $\partial M$ at
a right angle, a subarc $\gamma'$ of $\gamma$ would connect a truncation triangle
and its opposite hexagon in some tetrahedron $T'\in\calT$. Then the length
of $\gamma'$ would be greater than the distance between such a truncation triangle
and its opposite hexagon, which in turn is greater than the distance between $p$ and
some truncation triangle in $T$, since $T$ and $T'$ are isometric to each other and regular.
\end{proof}

\begin{rem}\label{alternative:rem}
{\em An alternative  proof of Theorem~\ref{distinct:teo}
could be based on the machinery developed in~\cite{FriPe}:
Remark~\ref{decomposition:rem} and the tilt-formula~\cite{Wee, Ush, FriPe}
easily imply that a one-edged ideal triangulation of a manifold
$M\in\calM_n$ is combinatorially equivalent to Kojima's canonical
decomposition of $M$, which is obtained by straightening the ideal triangulation 
dual to the cut-locus of $M$ (see~\cite{Ko}).}
\end{rem}

We say that a standard polyhedron is \emph{orientable} if it can be embedded
in an orientable 3-manifold.

\begin{cor}
The set $\calM_n$ is in one-to-one correspondence with the set of orientable standard polyhedra
with $n$ vertices and one face.
\end{cor}

Another remarkable consequence of Theorem~\ref{distinct:teo}
is that there is a finite (and very easy) algorithm
to decide whether an element of $\calM_n$
is chiral or amphichiral.  We defer the precise statement
to the next section (Proposition~\ref{chiral:prop}).

\section{Numerical estimates}

The results of the previous section would of course be of little or
no interest if $\calM_n$ (the class of manifolds having a triangulation
with one edge and $n$ tetrahedra) turned out to be empty or very small.
In this section we prove that $\#\calM_n$ grows at least exponentially with $n$, 
deducing that the number 
of orientable compact 3-manifolds of complexity $n$ also grows at least 
exponentially
(Corollary~\ref{exp:growth:cor}). We do this by concentrating on a special
class of one-edged triangulations, and we give some hints showing that our 
exponential lower estimate on $\#\calM_n$ is actually far from being sharp.
More accurate estimates would require a complicated combinatorial analysis
without providing a qualitatively better information on the rate of growth of
$\#\calM_n$.

\paragraph{Oriented spines and o-graphs} 
In the whole of this section we consider \emph{oriented} (rather than
just orientable) manifolds (but see Proposition~\ref{chiral:prop}).
We recall~\cite{BePe} that if $P$ is a standard spine of an oriented $M$
then $P$ also carries an orientation, defined as a screw-orientation
along the edges of $S(P)$ with a natural compatibility at vertices
(see~\cite[Fig.~2]{BePe:bis}). Conversely, if $P$ is an oriented standard polyhedron, then
$P$ is orientable, and the manifold it defines is oriented. In addition,
$P$ can be described by two additional structures on the 4-valent graph $S(P)$:
\begin{itemize}
\item an embedding in the plane of the neighbourhood of each vertex, 
with two opposite strands marked as being over the other two,
as in knot projections;
\item a colour in $\matZ/_3$ attached to each edge.
\end{itemize}
A 4-valent graph with these additional structures is called 
an \emph{o-graph}.
It was shown in~\cite{BePe} that any o-graph defines an oriented standard 
polyhedron, whence an oriented manifold, and that two o-graphs defining the 
same oriented polyhedron are related by certain ``C-moves.''
The effect of a C-move is
to change the planar structure 
at a vertex and the $\matZ/_3$-colouring of the edges incident to this vertex.

\begin{figure}
\begin{center}
\mettifig{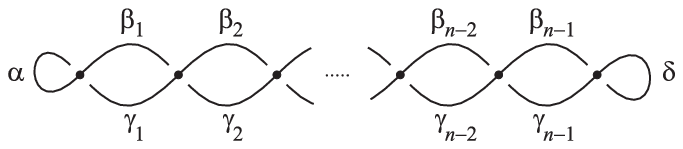,width=8cm}
\nota{O-graph of a generic polyhedron based on $G_n$}\label{open_chain_ogra:fig}
\end{center}
\end{figure}

\paragraph{O-graphs based on the open chain}
Let $G_n$ be the graph with vertices $v_1,\ldots,v_n$, a closed edge at $v_1$
and one at $v_n$, and two edges joining $v_{i}$ to $v_{i+1}$ for $i=1,\ldots,n-1$.
We characterize in this paragraph the oriented 
standard polyhedra $P$ such that $S(P)=G_n$ and
$P$ has a single face. We begin with the following fact that one can readily
establish using the C-moves mentioned above:

\begin{lemma} \label{open_chain:lem}
Any oriented standard polyhedron $P$ such that $S(P)=G_n$ can
be represented by an o-graph as shown in Fig.~\ref{open_chain_ogra:fig}.
\end{lemma}

The description of which o-graphs as in Fig.~\ref{open_chain_ogra:fig}
have a single face will use the language of finite state automata 
(see for instance~\cite{fsa-cite}).
We recall that a finite state automaton over a finite set $A$ (the alphabet)
consists of a finite set $S$ (the states), a function $S\times A\to S$ (the
transition function), an element $s_0$ of $S$ (the start state), and
a subset $S'$ of $S$ (the set of accept states). A word (a finite string of 
letters from the alphabet) is accepted by the automaton if, 
starting from $s_0$, reading the word
from left to right, and using the transition function,
the automaton ends in a state of $S'$. An automaton can be encoded by a picture,
where the states are represented by $S$-labeled boxes (with double margin
for accept states),
the transition function is given by box-to-box $A$-labeled arrows, and a 
mark indicates the start state.

\begin{teo}\label{open:chain:teo}
The oriented standard polyhedron defined by the o-graph of
Fig.~\ref{open_chain_ogra:fig} has a single face if and only if
both $\alpha$ and $\delta$ are different from $2$ and the word
$$(\beta_1,\gamma_1)(\beta_2,\gamma_2)\cdots(\beta_{n-1},\gamma_{n-1})$$
is accepted by the automaton described in Fig.~\ref{open_chain_fsa:fig}, where
\begin{figure}
\begin{center}
\mettifig{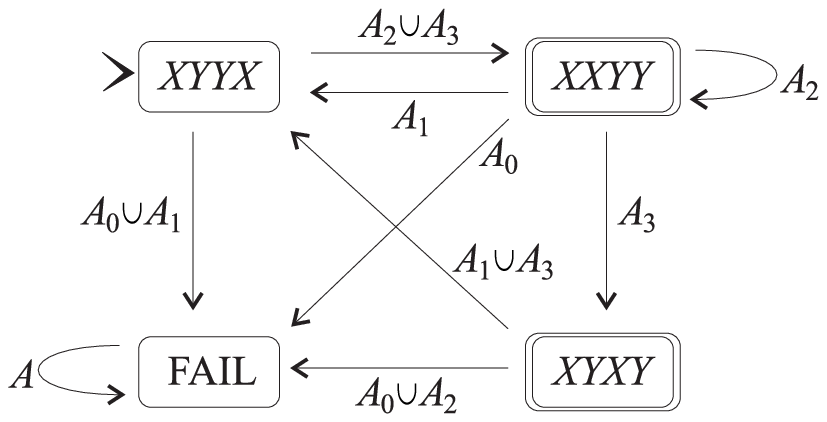,width=8cm}
\nota{A finite state automaton with alphabet $A=A_0\cup A_1\cup A_2\cup A_3$.} \label{open_chain_fsa:fig}
\end{center}
\end{figure}
$$
\begin{array}{l}
A_0=\{(2,2)\},\qquad A_1=\{(0,0),(1,1)\},\qquad
A_2=\{(1,0),(0,1)\},\\ 
A_3=\{(0,2),(1,2),(2,0),(2,1)\}, \qquad\quad\,
A=A_0\cup A_1\cup A_2\cup A_3.
\end{array}$$
\end{teo}

\begin{proof}
We confine ourselves to a general explanation, omitting the many combinatorial details.
We analyze the graph left to right, starting from $\alpha$.
Figure~\ref{dimo_chain_ogra:fig} shows that for $\alpha=2$ 
\begin{figure}
\begin{center}
\mettifig{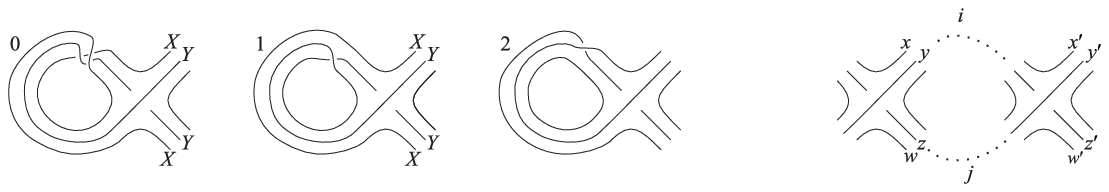,width=14cm}
\nota{Portions of o-graph.} \label{dimo_chain_ogra:fig}
\end{center}
\end{figure}
there are at least two faces, so $\alpha\in\{0,1\}$.
Now we examine the pair of colours $(\beta_1,\gamma_1)$ starting with 
the string $XYYX$ which describes the way the faces already constructed
are matched to the left of the
point we are considering. (Note that both $\alpha=0$ and $\alpha=1$ give $XYYX$.)
Depending on $(\beta_1,\gamma_1)$ we will have either the creation of a face
which does not fill the polyhedron, or a new pattern that describes the
matching of faces. 
More generally, as we proceed, we will have a string $xyzw$ of 
two symbols each repeated twice, and we will have to analyze the effect of a
pair of colours $(i,j)$, which either creates a closed face
or produces a new pattern $x'y'z'w'$ (see
Fig.~\ref{dimo_chain_ogra:fig} again). The detailed analysis of all the possibilities
leads precisely to the transitions shown in Fig.~\ref{open_chain_fsa:fig},
where ``FAIL'' means that a closed face is created. To conclude we note that
the final edge with colour $\delta\in\{0,1\}$ gives a single global face when
the input pattern is $XXYY$ or $XYXY$, and not otherwise.
\end{proof}

\paragraph{Growth of $\#\calM_n$}
Theorem~\ref{distinct:teo} shows that 
to compute $\#\calM_n$ it is sufficient to count the combinatorially distinct
standard polyhedra with $n$ vertices and one face. Restricting to the oriented
ones with the open chain $G_n$ as singular graph we must then discuss which o-graphs
as in Fig.~\ref{open_chain_ogra:fig} define the same polyhedron. 
A move that of course does not change the polyhedron associated to the
o-graph is the $180^\circ$ degree rotation.
Using the C-moves of~\cite{BePe} one can see that another such move
consists in interchanging each $\beta_k$ with the corresponding $\gamma_k$. 
In addition, these two moves are sufficient to generate
all graphs giving the same polyhedron. Therefore we have the following:

\begin{lemma}\label{only:four:lem}
An oriented standard polyhedron is defined by at most 
four different o-graphs as in Fig.~\ref{open_chain_ogra:fig}.
\end{lemma}

\begin{prop}\label{growth:estimate:prop}
There are at least $4\cdot 12^{(2n-5)/3}$ distinct oriented standard polyhedra
with one face and the open chain with $n$ vertices as singular set .
\end{prop}

\begin{proof}
We must count the possible choices for the $\beta_k$'s and $\gamma_k$'s,
\emph{i.e.}~the words of length $n-1$ accepted by the
finite state automaton of Fig.~\ref{open_chain_fsa:fig}, multiply by 4
(the choices for $\alpha$ and $\delta$), and divide by at most $4$ according
to the previous lemma. So it is sufficient to prove that there are at least
$4\cdot 12^{(2n-5)/3}$ words of length $n-1$ accepted by the automaton.

The idea is just to perform the loop $XYYX\to XXYY\to XYXY\to XYYX$ in all possible
ways, inserting a single loop $XXYY\to XXYY$ when $n-1$ is a multiple of $3$
(for in this case we would end up in the non-accept start state). Since
6 letters lead from $XYYX$ to $XXYY$, 4 lead from $XXYY$ to $XYXY$, and again $6$
from $XYXY$ to $XYYX$, it is clear that approximately 
$6^{2(n-1)/3}\cdot 4^{(n-1)/3}=12^{2(n-1)/3}$
words can be constructed with this method.
The exact computation carried out depending on the congruence class of $n-1$ 
modulo $3$ leads to the desired estimate.
\end{proof}

The next easy remark shows that the qualitative type of growth
just established is actually the maximal one could expect. 
After the remark we also give an obvious consequence of the 
previous proposition.

\begin{rem}
\emph{Given a 4-valent graph $G$ with $n$ vertices, there exist at most 
$2^n\cdot 3^{2n} = 18^n$
distinct oriented standard polyhedra $P$ such that $S(P)=G$. If $G=G_n$,
using Lemma~\ref{open_chain:lem}, one
can see that there exist at most $3^{2n}=9^n$ of them}.
\end{rem}

\begin{cor}\label{exp:growth:cor}
There exist $c>0$ and $b>1$ such that $\#\calM_n\geqslant c\cdot b^n$.
In particular, the number of distinct orientable,
boundary-irreducible, and acylindrical manifolds of complexity $n$
is at least $c\cdot b^n$.
\end{cor}

We remind the reader that the manifolds referred to in the previous
corollary are precisely those known to have standard minimal spines.
Now we have:

\begin{rem}
\emph{The number of distinct oriented standard polyhedra with $n$ vertices
is bounded from above by $18^n\cdot g(n)$, where $g(n)\leqslant (4n-1)!!$ is the
number of distinct four-valent graphs with $n$ vertices.}
\end{rem}

\paragraph{Further comments on estimates}
The lower bound on the number of elements of $\calM_n$ based on the open chain graph $G_n$
provided by Proposition~\ref{growth:estimate:prop}
is very far from being sharp. For instance, if we consider in Fig.~\ref{open_chain_fsa:fig}
the paths consisting of some $XYYX\to XXYY\to XYXY\to XYYX$ cycles
intermingled with some $XXYY\to XXYY$ loops, we deduce that the number of
distinct spines is at least
$$\sum_{k=0}^{n-2}\left\{
\begin{array}{ll}
0 & {\rm if}\ n-k-1=3h,\\
2^k\cdot 6^{2h+1} \cdot 4^h \cdot\left(\begin{array}{c} h+k \\ h \end{array}\right)
& {\rm if}\ n-k-1=3h+1,\\
2^k\cdot 6^{2h+1} \cdot 4^{h+1} \cdot\left(\begin{array}{c} h+k \\ h \end{array}\right)
& {\rm if}\ n-k-1=3h+2.\end{array}\right.$$
Concentrating on the term of the sum corresponding to $k=[n/7]$ and using
Stirling's formula one can for instance deduce from this estimate
that there exists $c>0$
such that, for all $\varepsilon>0$,
the number of oriented elements of $\calM_n$ based on $G_n$
is at least $c\cdot (6-\varepsilon)^n$ for $n\gg 0$.
Note that $12^{2/3}\cong 5.2418$.

\begin{rem}
\emph{With the aid of a computer, in~\cite{MaPe} we have
listed and classified the about 2,000
closed, irreducible and orientable manifolds of complexity up to $9$.
Corollary~\ref{exp:growth:cor} suggest that a similar listing for 
orientable, compact, boundary-irreducible, acylindrical 
manifolds may be hopeless. Our lower bound on their number, even if 
not sharp, already implies that there are at least
115,000 such manifolds in complexity up to 9.}
\end{rem}

\begin{figure}
\begin{center}
\mettifig{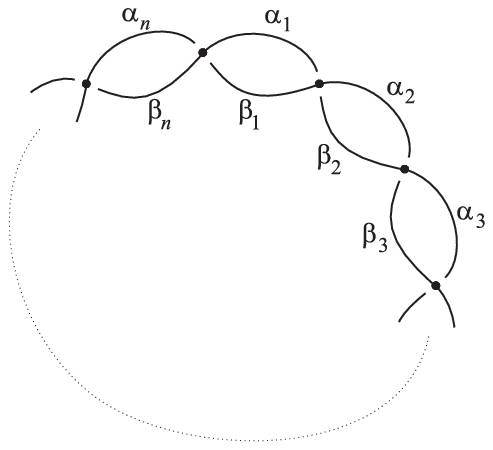,width=6.5cm}
\nota{An o-graph based on the closed chain.} \label{closed_chain_ogra:fig}
\end{center}
\end{figure}

\paragraph{Another example: the closed chain}
We have concentrated in this section on the open chain $G_n$, because
this graph is already sufficient to show that $\#\calM_n$ grows exponentially.
But we believe that a systematic investigation of the 4-valent graphs supporting 
a polyhedron with a single face would be quite interesting. Our guess is actually
that most graphs indeed support many different such polyhedra.
As another example we 
only mention the closed chain, shown in Fig.~\ref{closed_chain_ogra:fig}.
The combinatorial analysis is in this case harder than that carried out for
Theorem~\ref{open:chain:teo}, but we state at least the following fact, which already
implies that again in this case the number of relevant polyhedra grows 
exponentially with $n$.

\begin{prop}\label{closed_chain_prop}
Fix $(\alpha_1,\beta_1)\in\{(0,0),(0,1),(1,0),(1,1)\}$
in the o-graph of Fig.~\ref{closed_chain_ogra:fig}. Let 
$(\alpha_2,\beta_2)\cdots(\alpha_n,\beta_n)$ 
be a random word $w$ in the letters $\{(0,2),(1,2),(2,0),(2,1)\}$.
Then for $n\gg 0$ there is a probability $1/2$ that $w$
defines a polyhedron with a single face.
\end{prop}

\paragraph{Chirality} Using o-graphs, we discuss here chirality of
the elements of $\calM_n$, first in general and then for
those obtained from the open chain $G_n$.

Let $M\in\calM_n$ be oriented and let $\Gamma$ be an o-graph 
representing an oriented spine of $M$ with $n$ vertices.
According to~\cite{MaPe} the manifold $-M$ obtained by reversing
the orientation of $M$ is represented by the o-graph $-\Gamma$
obtained from $\Gamma$ by switching overstrands and understrands 
at vertices, and changing each edge-colour to its opposite in $\matZ/_3$.
Now $M$ and $-M$ have a unique oriented
spine with $n$ vertices by Theorem~\ref{distinct:teo}, so
$M$ and $-M$ are homeomorphic (\emph{i.e.}~$M$ is amphichiral) 
if and only if $\Gamma$ and $-\Gamma$
define isomorphic oriented polyhedra. Therefore:

\begin{prop}\label{chiral:prop}
$M$ is amphichiral if and only if $-\Gamma$ is combinatorially identical
to one of the $6^n$ o-graphs obtained from $\Gamma$ by C-moves at its vertices.
\end{prop}

Using C-moves one can now readily prove the following:

\begin{prop}
If $\Gamma$ is the o-graph of Fig.~\ref{open_chain_ogra:fig} then $-\Gamma$ is
a similar o-graph with colours:
$$\begin{array}{lll}
\alpha'=1-\alpha, & \beta'_{2k+1}=1-\gamma_{2k+1}, 
&\beta'_{2k}=1-\beta_{2k},\\
\delta'=1-\delta, & \gamma'_{2k+1}=1-\beta_{2k+1}, 
&\gamma'_{2k}=1-\gamma_{2k}.
\end{array}$$
\end{prop}

This proposition and Lemma~\ref{only:four:lem} show that
chirality for the elements of $\calM_n$ obtained from the open
chain $G_n$ can be tested very efficiently: we only need to check
whether one of four given $2n$-tuples of elements of 
$\matZ/_3$ coincides with another given such $2n$-tuple.

\vspace{1.5 cm}

\noindent
\hspace*{6cm}Scuola Normale Superiore\\ 
\hspace*{6cm}Piazza dei Cavalieri 7 \\
\hspace*{6cm}56127 Pisa, Italy\\ 
\hspace*{6cm}frigerio@sns.it
\vspace{.5 cm}

\noindent
\hspace*{6cm}Dipartimento di Matematica\\ 
\hspace*{6cm}Universit\`a di Pisa\\ 
\hspace*{6cm}Via F. Buonarroti 2\\ 
\hspace*{6cm}56127 Pisa, Italy\\ 
\hspace*{6cm}martelli@mail.dm.unipi.it
\vspace{.5 cm} 

\noindent
\hspace*{6cm}Dipartimento di Matematica Applicata\\ 
\hspace*{6cm}Universit\`a di Pisa\\
\hspace*{6cm}Via Bonanno Pisano 25B\\ 
\hspace*{6cm}56126 Pisa, Italy\\ 
\hspace*{6cm}petronio@dm.unipi.it

\end{document}

%% file: defs_oneface.tex
\newtheorem{lemma}{Lemma}[section] 
\newtheorem{teo}[lemma]{Theorem}
\newtheorem{rem}[lemma]{Remark} 
\newtheorem{prop}[lemma]{Proposition}
\newtheorem{cor}[lemma]{Corollary}

\newcommand{\matN}{\ensuremath {\mathbb{N}}}
\newcommand{\matR} {\ensuremath {\mathbb{R}}}

\newcommand{\matZ} {\ensuremath {\mathbb{Z}}}
\newcommand{\matC} {\ensuremath {\mathbb{C}}}

\newcommand{\matH} {\ensuremath {\mathbb{H}}}

\newcommand{\calM} {\ensuremath {\mathcal{M}}}

\newcommand{\calT} {\ensuremath {\mathcal{T}}}

\newcommand{\TV}{{\rm TV}}

\newcommand{\nota} [1] {\caption{\footnotesize{#1}}}

\newfont{\Got}{eufm10 scaled 1200}

\font\titsc=cmcsc10 scaled 1200

\newcommand{\dimo}[1]{\vspace{2pt}\noindent\textit{Proof of \ref{#1}}.\ }
\newcommand{\finedimo}{{\hfill\hbox{$\square$}\vspace{2pt}}}